\documentclass[14pt]{article}
\usepackage[english]{babel}
\usepackage{amsthm,amsfonts, amsbsy, amssymb,amsmath,graphicx}
\usepackage{graphics}
\author{Vassily Olegovich Manturov
\footnote{Bauman Moscow State Technical University and Chelyabinsk State University}  
}

\date{}
 \theoremstyle{plain}
\newtheorem{thm}{Theorem}

\theoremstyle{dfn}
\newtheorem{dfn}{Definition}
\newtheorem{rk}{Remark}
\newtheorem{example}{Example}

\newtheorem{st}{Statement}

 \def\Z{{\mathbb Z}}
 \def\0{{\mathbbf 0}}
 \def\1{{\mathbbf 1}}

 \def\gG{{\mathfrak G}}
 \newcommand{\skcr}{\raisebox{-0.25\height}{\includegraphics[width=0.5cm]{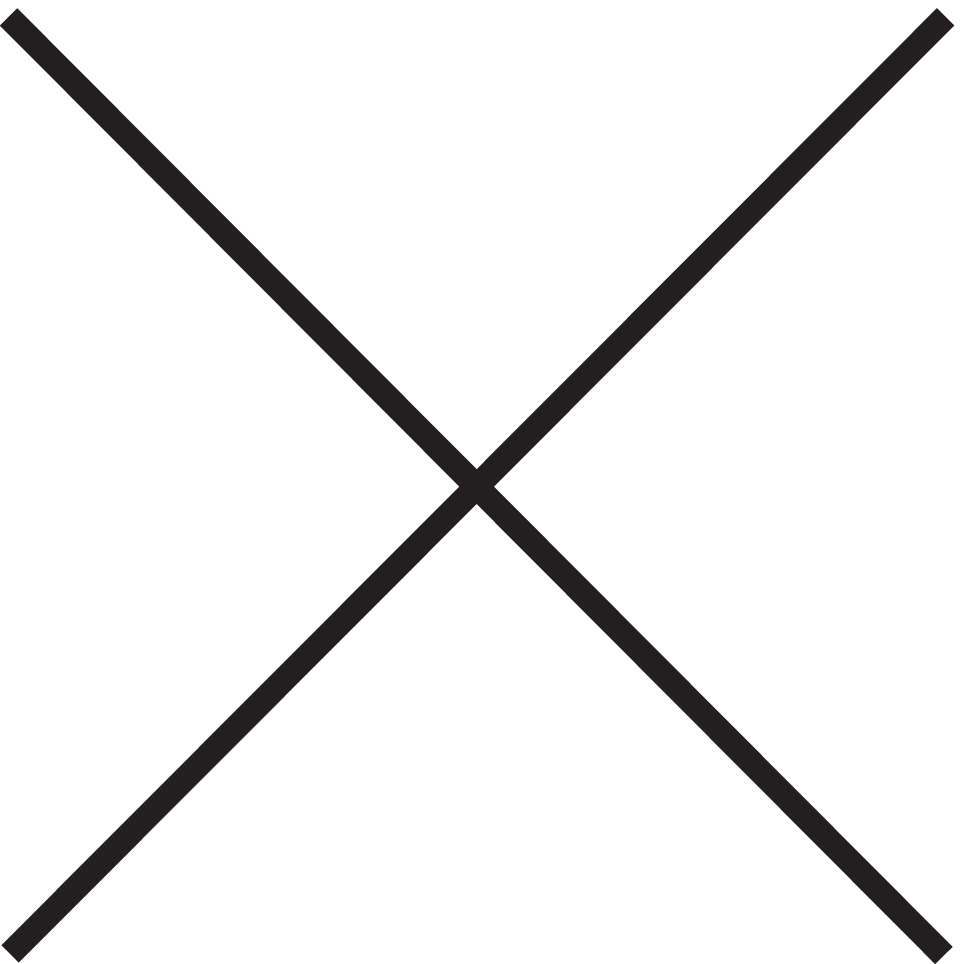}}}
 \newcommand{\skcrv}{\raisebox{-0.25\height}{\includegraphics[width=0.5cm]{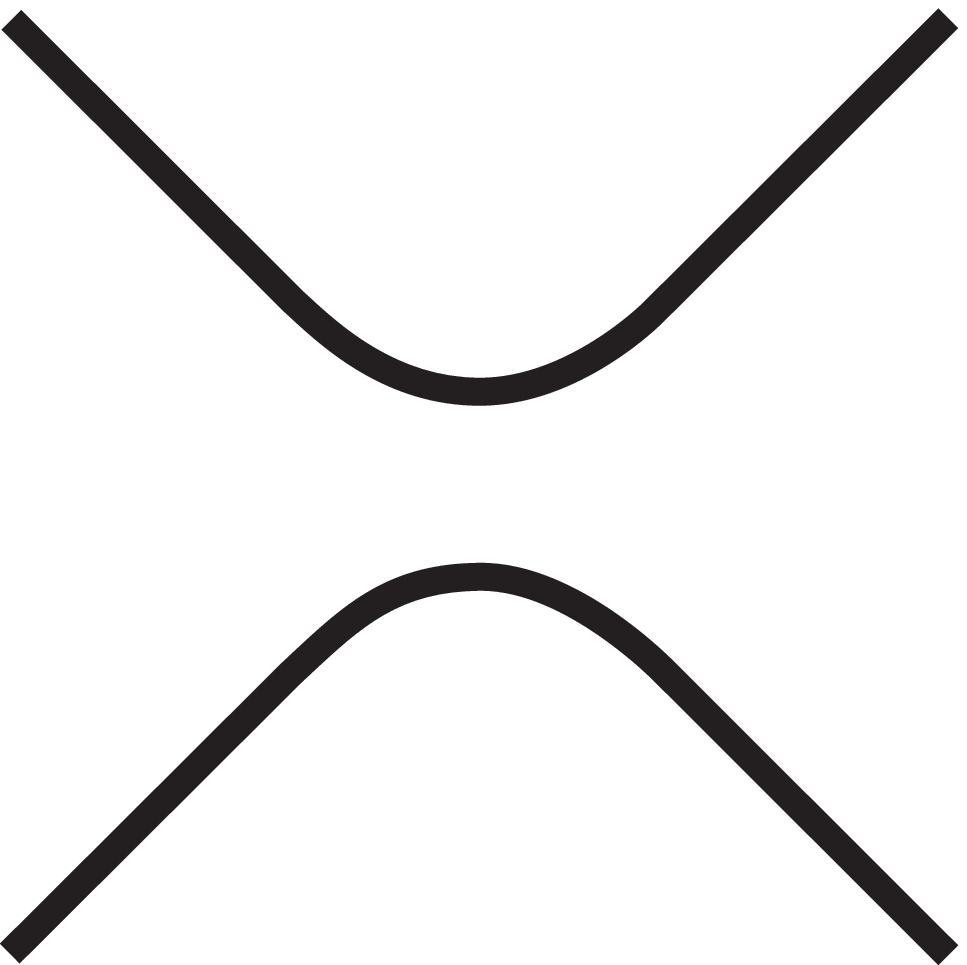}}}
\newcommand{\skcrh}{\raisebox{-0.25\height}{\includegraphics[width=0.5cm]{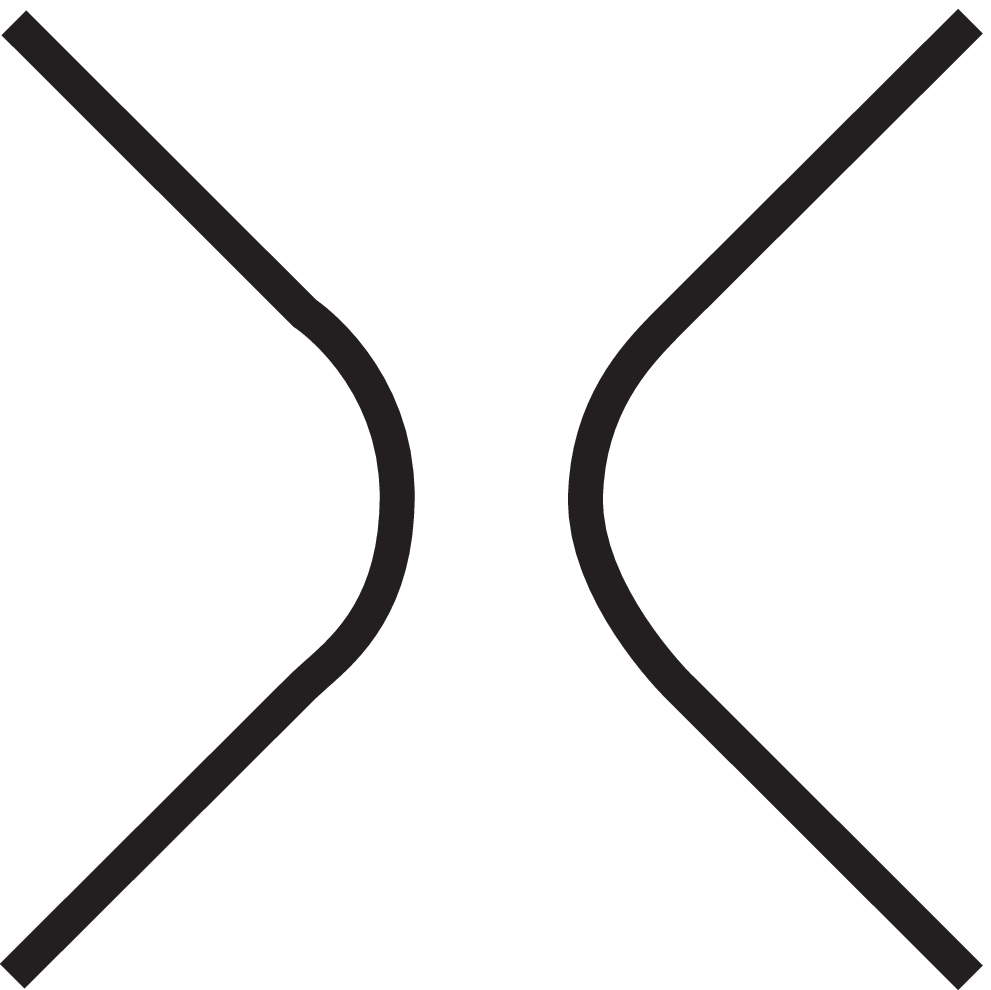}}}
\newcommand{\skcrosso}{\raisebox{-0.25\height}{\includegraphics[width=0.5cm]{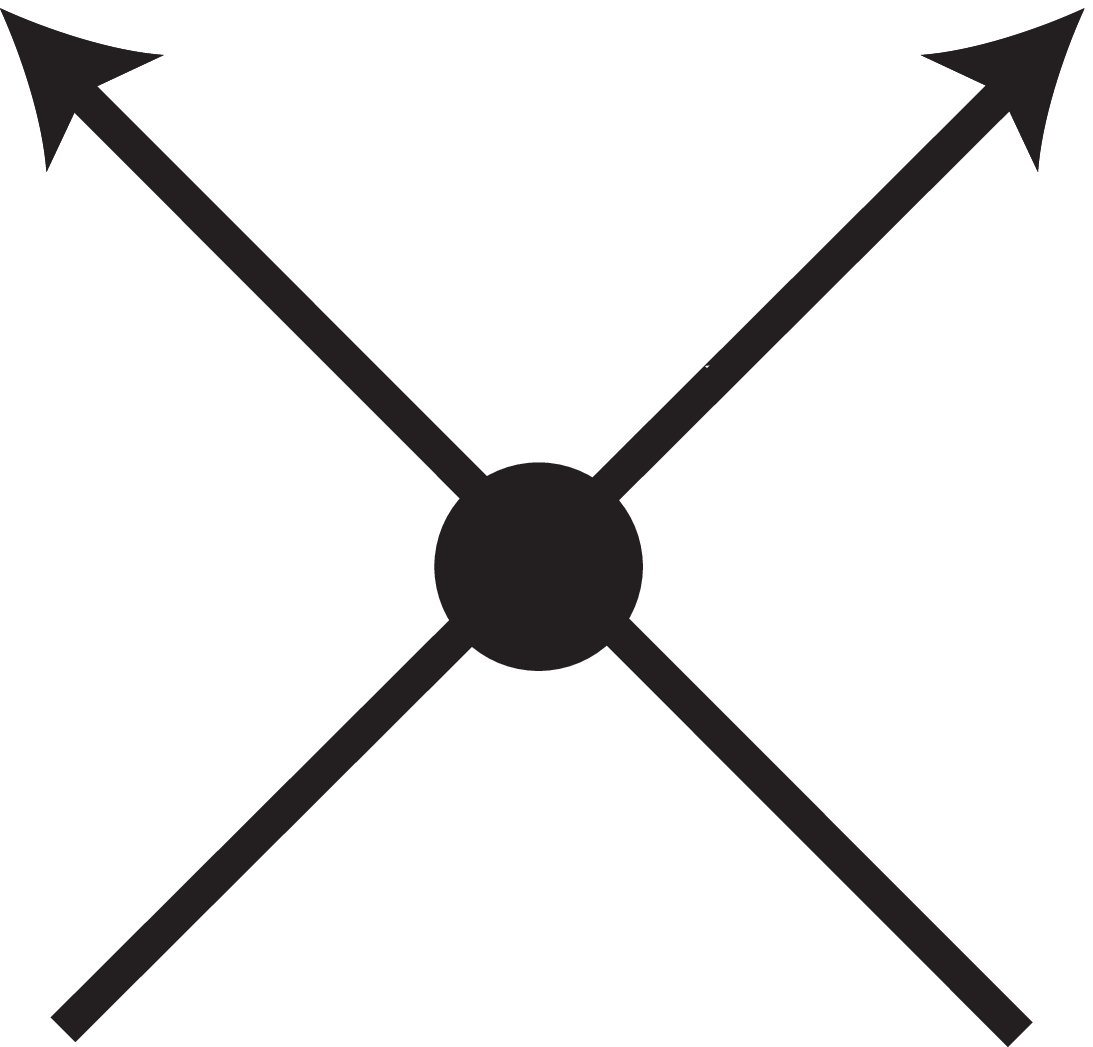}}}

 \newcommand{\skcrho}{\raisebox{-0.25\height}{\includegraphics[width=0.5cm]{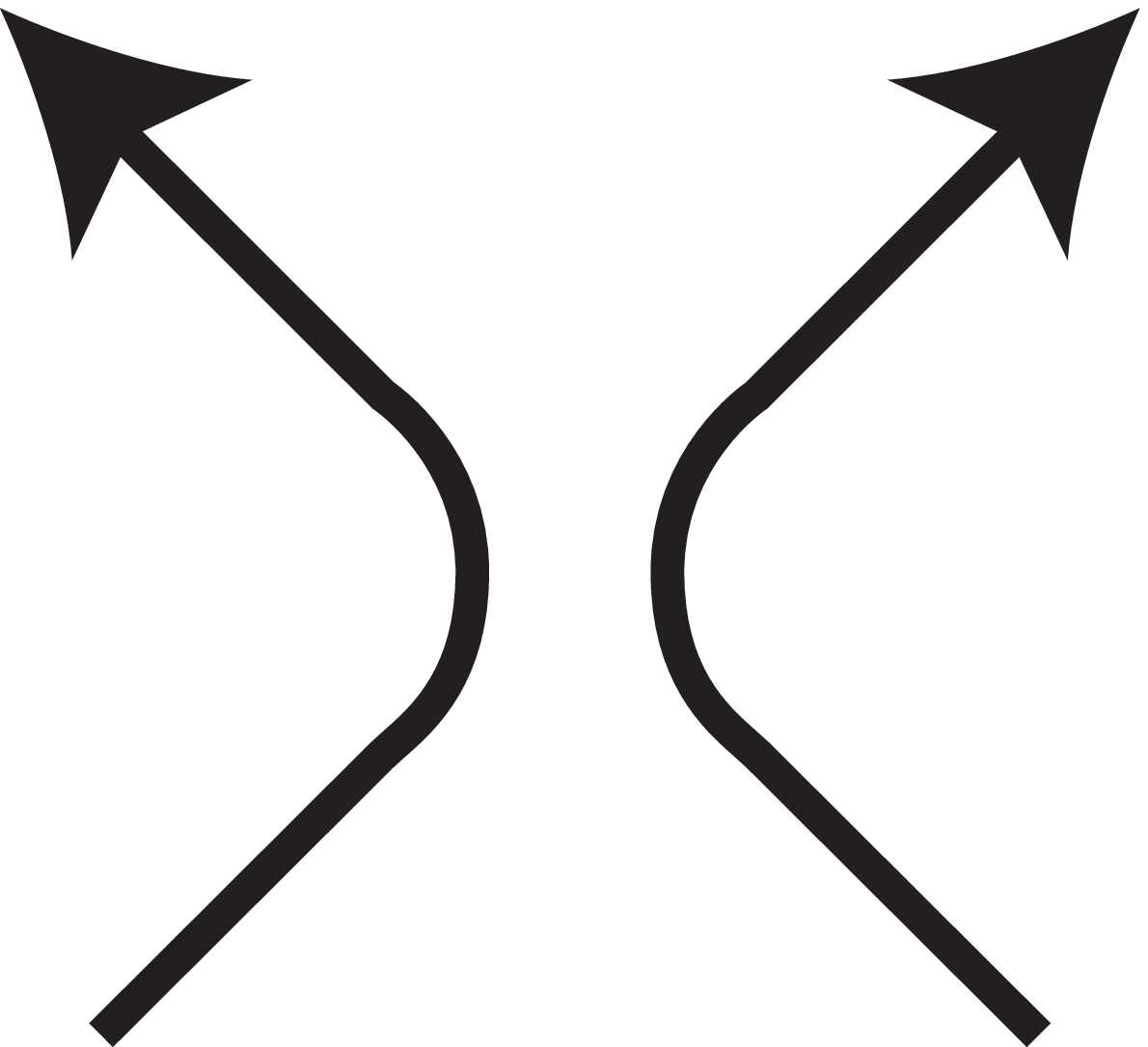}}}

\title{New Parities and Coverings over Free Knots}
\date{}

\begin{document}

\maketitle

\begin{center}

{\bf Abstract}
\vspace{0.5cm}

%\begin{abstract}

{\large{
In the present paper, we develop the parity theory invented in \cite{ManSb};
we construct new parities for two-component (virtual and free) links. New parities significantly
depend on geometrical properties of diagrams; in particular, they are mutation-sensitive.

New parities can be used practically in all problems, where parities were previously
applied.
}
}
\end{center}
%\end{abstract}


\section{Introduction}

The theory of multi-component links is much richer than the theory of knots.
Thus, for classical links there exists a powerful link-homotopy theory
\cite{Milnor,HL}: the simplest link homotopy invariant is the linking coefficient. Classical
invariants of knots can be modified if one passes from knots to multicomponent links.

Virtual knot theory invented by Kauffman \cite{KaV} includes classical knot theory.
As it turns out, the behaviour of one-component virtual knots is similar to the behaviour of classical {\em links};
virtual knots admit invariants similar to invariants of multicomponent classical links. Thus,
in \cite{KaSelfLink} the ``self-linking coefficient'' for
virtual knots was defined.

This effect is possible because of {\em parities} and {\em coverings} over virtual knots: there exists a
well defined map from the set of virtual knots to the set of virtual links.

The parity theory invented by the author \cite{ManSb}, see also
\cite{IMNObzor,ManSb2}, allows one to make a distinction between
two types of crossings of virtual knots, {\em even ones} and {\em odd ones} (in fact, this can be done
even for a simplification of virtual knots called free knots); this circumstance allows one to treat
{\em virtual knots} as {\em multicomponent links} by means of coverings coming from parities. These coverings
were first used in \cite{ManKhNakr} for the construction of Khovanov homology theory for virtual knots.

Usually, knot invariants are valued in some algebraic objects: numbers, polynomials,
groups. The information about knots encoded by such invariants is very implicit. Say,
a grading or a degree of some knot invariant allows one to judge about the crossing number
of a knot, but not about the geometrical form of diagrams of the knot.
The algebraic nature of invariants allows one to construct easily some
transformation which do not change the values of the invariant,
but change the knot dramatically.

For an important example one can take {\em mutations}, \cite{CrM}. Mutations change the topology
of the knot dramatically but do not affect the value of most of algebraic invariants.

The parity theory allows one to construct knot invariants valued in {\em knot diagrams}, which allows one to solve problems of the following sorts:
\begin{enumerate}
\item
Reduce properties of knots to properties of their diagrams;

\item Construct functorial mappings between knots.

\end{enumerate}

Particularily, the first problem can be solved by using the so-called ``parity bracket'' introduced by the author in
\cite{ManSb}. This parity bracket is a well defined map from the set of knots to the set of graphs;
if some ``oddness'' and ``irreducibility'' conditions are satisfied, then the bracket of a free knot diagram is equal to the diagram itself.
%%%%%%%%%%%%%%%%%

The first statement can be formulated as the following ``meta''theorem

\begin{thm} If a (virtual) knot diagram $K$ is complicated enough (``odd'') then any diagram $K'$ equivalent to it contains $K$ as a subdiagram.\label{meta}
\end{thm}

\begin{figure}
\centering\includegraphics[width=200pt]{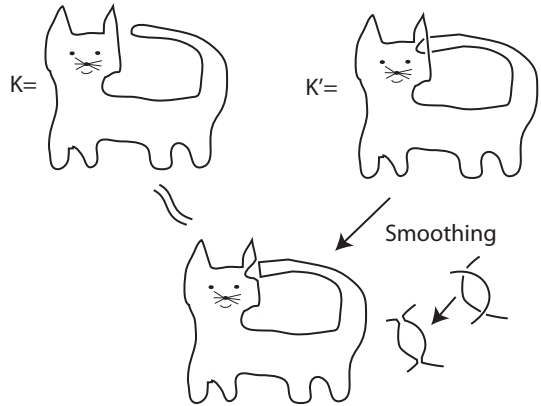}
\caption{Diagram $K$ is a subdiagram of some diagram $K'$ equivalent to $K$} \label{cat}
\end{figure}

View Fig. \ref{cat}. This figure demonstrates how the theorem applies to a unique second Reidemeister move
applied to the diagram $K$.

Among parities, the {\em Gaussian parity} plays a special role. It is the only non-trivial parity for free knots. This parity is also defined for knotted spheres (see
\cite{ManSb2}), which are smoothly mapped to $3$-manifolds. A formal definition in \cite{ManSb2} was given for discs,
 but it can be generalized verbatim for the case of spheres. Thus, the Gaussian parity allows one to construct not only invariants of free knots but also {\em a sliceness obstruction}: an obstruction to span a free knot by a disc with standard ``$3$-dimensional'' singularities.

 The other parities for knot theories require some additional structure defined on knots/links; there is no such general structure for free knots.

The major part of invariants constructed by the author in previous papers (for references, see
\cite{MI,IMNObzor}) are well defined for any parity.

However, the parities themselves for knot theories can be constructed easily by looking at some combinatorial or homological properties. Parities that we are going to construct in the present paper are based on some ``patterns'' (subdiagrams) $P$ in a given diagram. This means that some crossing of a knot diagram is decreed even or odd depending on whether there exist some ``subknots'' inside the knot located in a predicted manner with respect to the given crossing. This makes parity arguments even more delicate.

\begin{figure}
\centering\includegraphics[width=3.0in]{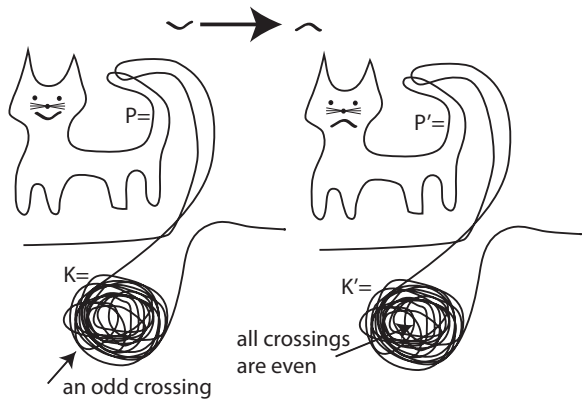}
\caption{A parity coming from a pattern} \label{cat2}
\end{figure}

In  Fig. \ref{cat2}, we schematically show a way of changing a link or its component (say, mutation); together with its pattern in question $P\to P'$;  such a transformation will not change ``algebraic'' invariants whose construction is based on some ``bare count``; however, such an effect can be tackled by looking at parities based on patterns.

The new parities are based on some count of intersections with a {\em given pattern $P$}.

Thus, if the sample changes slightly as shown in Fig.~\ref{cat2},
then the parity changes drastically. For example, the right part of
Fig.~\ref{cat2} does not contain $P$ at all, so all the crossings
are even, whereas some crossings in the left part of the figure are
odd.

Let us note that various topological constructions of parity theory
lead us to ``graphical'' results in classical knot theory,
see~\cite{CM,KM,Ch}.

In the present paper, we shall show how to construct various new parities for the two-fold covering over a free knot, where the two-fold covering is constructed by using the Gaussian parity, and how to use these new parities to construct new invariants of free knots.

In the sequel, we are planning to use new parities in order to construct invariants of free knot cobordisms.

The present paper is organized as follows.

In the second section, we give the definition of free knots, known parities, the parity bracket, and present a detailed discussion of the Turaev delta (also known as Turaev's bracket) as well as a detailed discussion of the two-fold covering over virtual knots.

In the third section, we define the parity for crossings of one component in a two-component link and show that in the case of a two-component link, the bracket and the delta may lead one to more elaborated parities.

The last section is devoted to the construction of various examples illustrating various subtleties of new parities and invariants coming from them.

\subsection{Acknowledgements}

The present work is supported by the Laboratory of Quantum Topology ot the Chelyabinsk State University (Russian Govenrment Grant No,
; 14.Z50.31.0020)

I am deeply indebted to I.M.Nikonov for various discussions of various (Russian and English) versions of
the text; I am also grateful to D.P.Ilyutko, D.A.Fedoseev, and V.A.Krasnov for fruitful discussions.

I am also grateful to the referee for various remarks concerning the text.

\section{Basic definitions}

\subsection{Free knots and links}

\begin{dfn}
By a
{\em $4$-graph} we mean a $1$-dimensional complex each component of which being either homeomorphic to the circle or being a four-valent graph; by {\em vertices} of a $4$-graph we mean vertices of its components which are not homeomorphic to circles; by {\em edges} of a $4$-graph we mean edges of its components which are not homeomorphic to circles as well as components homeomorphic to the cirlce (the latter will be called {\em cyclic components} or {\em cyclic edges}); every non-cyclic edge will be treated as an equivalence class of the two (different!) half-edges composing it.
\end{dfn}

\begin{rk}
Note that the CW-structure is immaterial in the sequel; we assume that there are no vertices on cyclic edges.
\end{rk}

\begin{dfn}  We say that a  $4$-graph is  {\em framed}, if at each
vertex of this graph, the four half-edges are split into two pairs of {\em opposite} ones;
we shall call such a splitting a {\em framing}; half-edges which are incident to the same vertex and are not opposite, will be called {\em adjacent}. Framed $4$-graphs are considered up to the natural {\em equivalence}, i.e., a homeomorphism which preserves the framing (we do not impose any restrictions on cyclic components).
\end{dfn}

\begin{rk} In the sequel, we treat framed $4$-graphs only up to equivalence.
\end{rk}

By a {\em cyclic universal component} of a framed $4$-graph we mean an equivalence class of edges of this graph containing
only one cyclic edge.

\begin{figure}
\centering\includegraphics[width=100pt]{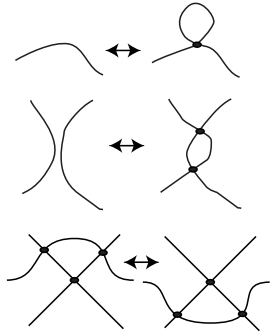}
\caption{Reidemeister moves} \label{rmoves}
\end{figure}

\begin{rk}
When drawing framed $4$-graphs on the plane, we shall not
indicate the framing assuming that formally opposite half-edges are exactly those locally opposite
on the plane. Framed $4$-graphs appear naturally when considering generic immersions
of a collection of circles into the plane; the framing is inherited from the immersion.

Homotopy classes of curves in $2$-surfaces motivate the following set of elementary equivalences on framed $4$-graphs called {\em Reidemeister moves}, see Fig. \ref{rmoves}.

When drawing diagrams on the plane, we will show only their changing parts (e.g., for moves).
\end{rk}

\begin{dfn} By a {\em free link} \cite{ManSb} we mean an equivalence class of framed $4$-graphs modulo Reidemeister moves.
\end{dfn}

\begin{dfn}
By an
{\em orientation} of a free link we mean orientations of all its cyclic components together with orientation of all edges belonging to non-cyclic components in such a way that for each two edges which are incident to the same vertex and opposite ate this vertex, one is incoming and the other one is emanating.
\end{dfn}

We say that a Reidemeister move is {\em increasing}, if it increases the number of vertices and {\em decreasing} if it decreases the number of vertices.

By a
 {\em non-cyclic unicursal component} of a framed  $4$-graph we mean an equivalence class of its edges generated by the elementary equivalence: two edges are elementary equivalent if they are opposite at some vertex.

One naturally gets  (see \cite{ManSb}) the
{\em number of unicursal components} of a free link.

Each crossing of a free link belongs either to one unicursal component or to two different components. In the former case, we shall call it {\em the pure crossing}, in the latter case, the crossing is {\em mixed}.

\begin{figure}
\centering\includegraphics[width=200pt]{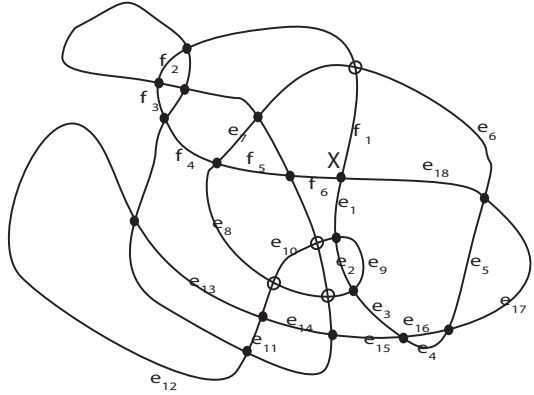}
\caption{A $2$-component free link and halves corresponding to a vertex} \label{halves}
\end{figure}

If a crossing $X$ of a framed  $4$-graph $\Gamma$ is not mixed then the unicursal component
passing through this crossing is naturally split into two {\em halves}. We say that the two edges  $a,b$ belonging to one unicursal component belong to the same {\em half} if there exists a chain of edges $a=e_{1},e_{2},\cdots, e_{k}=b$ where every two adjacent edges $e_{l},e_{l+1}$ are opposite at some vertex distinct from $x$. In Fig. \ref{halves}, chords
of the chord diagram corresponding to one half with respect to the vertex $X$, are denoted by $e_{i}$, and chords of the other half are denoted by $f_{i}$.

We say that a ($2$-component) link is {\em splitting} if it can be represented by a splitting diagram;
a ($2$-component) diagram, in turn, is {\em splitting} if it has no mixed classical crossings.

\begin{figure}
\centering\includegraphics[width=200pt]{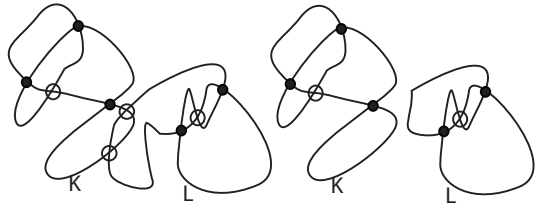}
\caption{Splitting diagrams} \label{split}
\end{figure}

In Fig. \ref{split}, both diagrams are splitting ones. The diagram on the LHS is splitting because it has no pure classical crossings; the diagram on the RHS has no mixed crossings at all; those ``virtual crossings'' depicted by circles are not counted as vertices.

\begin{rk}
Generically,
diagrams of curves on two-surfaces, carry more information than framed $4$-graph. To a framed $4$-graph,
there corresponds many different immersions into surfaces  (even if we restrict ourselves to connected  $4$-graphs
and even if we require that the graph tiles the surface into cells)). Thus, there is a natural ``forgetful'' map from
homotopy classes of classes of curves on $2$-surfacest to framed $4$-graph.

This map exactly corresponds to the forgetful map from flat links to free links.
\end{rk}

\subsection{The Parity. The Gaussian Parity}

Assume some class of knots $\mathcal{K}$ represented by equivalence classes of diagrams modulo moves is given;
the diagrams are required to be framed $4$-graphs, possibly, with some additional structures and/or restrictions, and
possibly, some additional restrictions are imposed on Reidemeister moves. Thus, classical knots can be thought of as {\em planar} diagrams modulo moves where {\em edge cyclic structure} together with {\em ovepass/underpass} structure is indicated, and natural restrictions are imposed on the second and the third Reidemeister moves.

Assume we are given a  {\em rule} which assigns with each vertex of a diagram $K$ from the class  $\mathcal{K}$ a  number $0$ (such vertices will be called {\em even})\index{Vertex!even} or $1$ (in this case, the vertex is called  {\em odd}).\index{Vertex!idd} In the sequel, we shall denote the parity of a vertex $v$ by $p(v)$ and write
 $p(v)=0$ if the crossing $v$ is even and $p(v)=1$ if the crossing $v$ is odd.

 \begin{dfn}
 We say that the above rule satisfies the parity axiomatics if the following holds. For every two diagrams  $K_1$ and $K_2$
 obtained from each other by one Reidemeister moves where  $K_2$ has no more crossings than $K_1$ the following conditions hold:
  \begin{enumerate}
   \item
if $K_2$ is obtained from $K_1$ by the first Reidemeister move, then the crossing of $K_1$ participating in the Reidemeister move, is {\em even};
   \item
   if  $K_2$ is obtained from $K_1$ by the second Reidemeister move, then both crossings taking part in the Reidemeister move {\em have the same parity};
   \item
   if $K_2$ is obtained from $K_1$ by the third Reidemeister move, then there exists a natural correspondence between the triple of crossings of $K_1$ and the triple of crossings of $K_2$ taking part in the Reidemeister moves($(a,a'),(b,b'),(c,c')$), see Fig.~\ref{sootvt}.

 \begin{figure}
\centering\includegraphics[width=200pt]{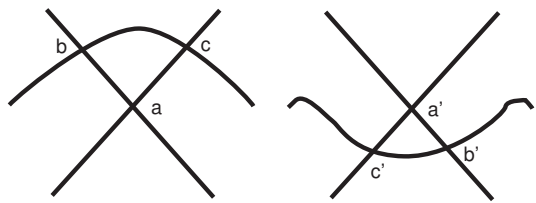}
\caption{The third Reidemeister moves and corresponding crossings} \label{sootvt}
 \end{figure}

We require that

a)
$p(a)=p(a'),p(b)=p(b'),p(c)=p(c')$,

b) among $a,\,b,\,c$, the number of odd crossings is even (i.e., equals \ $0$ or $2$).
  \item
We require that those crossings which remain untouched when passing from $K_{1}\mapsto K_{2}$, preserve their parity.
  \end{enumerate}\label{parityDef}
 \end{dfn}

 \begin{rk}
 If a framed
$4$-graph admits a symmetry (i.e., a framing-preserving isomorphism), then it follows from the definition that this symmetry preserves the crossing parity.
 \end{rk}

As an example, let us define the {\em Gaussian parity}, as follows.

 \begin {dfn}
A {\em chord diagram}\index{Diagram!chord} is a regular $3$-graph consisting of a selected non-oriented Hamiltonian cycle (the {\em core}) of the chord diagram \index{Core} and non-oriented edges ({\em
chords}),\index{Chord}, where each chord connects some two points on the circle, in such a way that different chords do not share points on the circle.

We say that two chords are {\em linked},\index{Chords!linked} if the endpoints of one of them belong to different connected components of the set obtained from $S^1$ by deleting the endpoints of the other chords. Otherwise, the chords are called {\em unlinked}.

We shall also admit the {\em empty chord diagram}, i.e., the chord diagram without chords which is formally not a graph..
 \end {dfn}

\begin{rk}
Analogously, one considers {\em multichord diagrams} or {\em chord diagrams on several circles}.
\end{rk}

\begin{dfn} A chord of a chord diagram is {\em even} (with respect to Gauss) is the number of chords, it is linked with, is even. Otherwise, we call this chord {\em odd}. A vertex of a framed $4$-graph is {\em even} or {\em odd} depending on whether the corresponding chord of the chord diagram is even or odd.

Finally, we say that a chord diagram is {\em even} if all chords of it are {\em even} and {\em odd} if all chords of it are odd.

\end{dfn}

If a framed $4$-graph $\Gamma$ has one unicursal component, then it can be thought of as the image of the map $f:S^{1}\to \Gamma$, which is bijective everywhere outside neighbourhooods of vertices of $\Gamma$ and their preimages.

Thus, framed $4$-graphs we are interested in are encoded by chord diagrams, and chords correspond to vertices.

\begin{figure}
\centering\includegraphics[width=200pt]{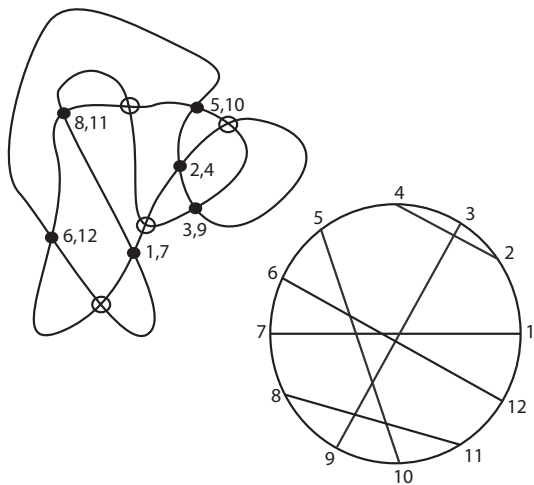}
\caption{A framed $4$-graph and the corresponding chord diagram}
\label{graphcd}
\end{figure}

In Fig. \ref{graphcd}, we show a framed $4$-graph and the chord diagram corresponding to it.

Absolutely analogously, one can define the correspondence between framed
$4$-graphs on many components and multichord diagrams (chord diagrams on many circles).

If we deal with multicomponent chord diagrams, parity can be defined only for pure crossings.

\begin{dfn} We say that a vertex of framed chord diagram is {\em even} (with respect to Gauss) if the corresponding chord of the chord diagram is even, and {\em odd}, otherwise.\end{dfn}

In other words, one can define the parity of a vertex $X$ (we deal with a free knot) of a framed $4$-graph $G$ as the parity of the number of crossings (common vertices) between the two halves of $G$ formed by $X$.

It can be easily seen that for free knots, the Gaussian parity satisfies all parity axioms.

In  \cite{IMNObzor} it is proved that the Gaussian parity is the only parity for free knot.
Parities for free links are also classified there.

The aim of the present paper is to find non-trivial ``parities'' which correspond to pure crossings of a free link only, provided that a link has a second component. From the formal point of view, this parity is not a parity for two-component links since it is defined {\em not for all crossings}; however, as we shall see further, there will be infinitely many parities of such sort, moreover, non-trivial free knot homotopy type will lead to such parities by means of ``patterns'', and these parities will depend on ``geometry'' of both components.

\subsection{The Bracket and The Delta}

In the present section we recall the two important mapping for free knots, {\em the parity bracket}, see \cite{ManSb}, and {\em the Turaev Delta}, see Fig. \cite{Tur}.

The first one gives an invariant of free knots valued in linear combinations of framed $4$-graphs.

The second one is a mapping from free knots (with one component) to free links (with two components).

For each framed $4$-graph $\Gamma$ with a non-empty set of vertices, for each vertex $x$, we define the two {\em splittings}
$\skcr\to \skcrv$,$\skcr\to \skcrh$, both being framed $4$-graphs.

An {\em even} splitting is a {\em splitting} performed at one or several even crossings.

Let us define the linear space $\gG$ as the set of
$\mathbb{Z}_{2}$-linear combinations of the following objects. We consider all framed $4$-graphs modulo the following equivalence relations:

1) the second Reidemeister move;

2) $L\sqcup \bigcirc=0$, i.e., the framed $4$-graph having more than one component with at least one trivial component, is assumed to be equal to zero.

Denote by $\gG'$ the linear subspace of the space $\gG$ generated by framed $4$-graphs with a unique unicursal component.

\begin{rk}
One can easily see that the equivalence classes listed below are uniquely characterized by their minimal representatives, i.e., framed  $4$-graphs without single circles and bigons. Thus, for summands from $\gG$ and $\gG'$, we shall use the term ``graph''
assuming some minimal element of the corresponding equivalence class.
\end{rk}

There exists a natural map $g\colon\gG\to \gG'$, taking all equivalence classes of framed graphs with more than one unicursal component, to zero. It is evident that $g$ is a group epimorphism.

Let $K$ be a framed-graph such that each unicursal component of $K$ has evenly many edges.
In particular, the construction given below will work for the case of free knots: if we have exactly one component then all edges belong to it.

 The invariant $\{\cdot\}$ (see
\cite{ManSb}) is given by
 $$
\{K\}=\sum_{s_{{even}}}K_{s_{even}}\in \gG,
 $$
 where the sum is taken over all even smoothings $s_{even}$ of the framed $4$-graph $K$,
 which are considered as elements of the group $\gG$.

 Here by an {\em even} smoothing we mean a smoothing at all {\em even} crossings.

 \begin{thm}
 The bracket $\{\cdot\}$ is an invariant of free links.\label{cdot}
 \end{thm}

The bracket $[\cdot]$ is defined for free knots and is obtained by projecting the bracket $\{\cdot\}$ to the subspace $\gG'$; in other words,
$$[K]=g(\{K\}).$$

Let us now modify the bracket $[\cdot]$ as follows. We shall consider the two-component free links $K\cup L$ and apply the smoothings to even crossings of the second component, i.e., we set
$$
[K\cup L]_{L}=\sum_{s_{{even}}}K\cup (L_{s_{even, 1 comp}}),
 $$
 where the sum is taken over all even smoothings $s_{even}$ of the framed
$4$-graph $L$ having one component.

We shall get a collection of graphs having two components each: one component is exactly the component  $K$, and the other component is $L_{s}$ depending on the state $s$. By looking at components $L_{s}$, one can construct further parities of the component  $K$.

\subsubsection{The map (Turaev's cobracket~\cite {Tur})}

We shall construct a map from $\Z_{2}$-linear combinations of oriented free knots to  $\Z_{2}$-linear combinations of oriented non-split free links.

Given an oriented framed graph $\Gamma$, and let
$X_{1},\dots, X_{N}$ be the crossings of $\Gamma$. At each crossing $X_{i}$ there is a unique way of {\em smoothing} for $\Gamma$ respecting the orientation: $\skcrosso\to\skcrho$; we shall denote the result of smoothing by $\Gamma_{i}$. Thus, we can define a map
 \begin{equation}
\Gamma\mapsto \sum_{i}\Gamma_{i},\label{equaa}
 \end{equation}
 where the sum is taken over those crossings for which the link
 $\Gamma_{i}$ is non-split.

The restriction above on $\Gamma_{i}$ is imposed for the map to be well defined: if we apply the first increasing Reidemeister move to $\Gamma$, then on the RHS of (\ref{equaa}) we shall get split links with trivial components as additional summands.

Thus, in order to get a well defined delta (cobracket), one should weaken some conditions, and forbid not all split summands, but only those where one component is trivial.

The invariance of delta under the second and the third moves is checked straightforwardly, see \cite{ManSb}.

Among the summands (\ref{equaa}), one can naturally take only those where the resulting links looks very specially (e.g., one component represents some concrete free knot).

Analogously, one can define maps from $n$-component free links to $\Z_{2}$-linear combinations of
$(n+1)$-component free links, where the corresponding sum is taken only over those crossings lying in one component, with some additional restrictions and factorizations.

In the sequel, we shall need various modifications of the map
$\Delta$ to get free links and free knots and to construct parities.

\subsection{The projection and the twofold covering}

Each parity in a knot (link) theory induces two natural maps
(projection and a twofold covering) that are constructing in the
following way.

Let $K$ be a framed $4$-graph. We construct graphs $K^{2}$ and $K'$
as follows. Firstly, assume that the graph $K$ is connected.

If $K$ is a cycle, then the graph $K'$ consists of one cycle and the
graph $K^{2}$ consists of two cycles of the same length. We say that
the cycles of $K^{2}$ are {\em dual}\index{Cycles!dual} to each
other and that they {\em cover} the cycle of $K$.

To each vertex  $v$ of the graph $K$ we assign two {\em covering}
vertices $v_{1},\,v_{2}$.\index{Vertex!covering} These vertices will
be called {\em dual}.\index{Vertices!dual} Choose a spanning tree
$T$ of the graph $K$.

We call all the edges of $T$ to be {\em good}.\index{Edge!good}
Edges of $K$ which do not belong to $T$ are divided into {\em good}
and {\em bad} edges as follows.\index{Edge!bad} Any edge $e$ in
$G\setminus T$ connects vertices $v,\,w$ and determines in $K$ a
minimal cycle which consists of the edge $e$ and the shortest path
from $v$ to $w$ in the tree. This cycle is {\em
rotating}\index{Cycle!rotating} (from an incident edge to a
neighbouring edge) at some vertices and {\em
transversal}\index{Cycle!transversal} at the remaining vertices. We
call the edge $e$ {\em good} if the number of transversal vertices
in the cycle is even.

To each edge $e$ of the spanning tree $T$, which connects two
vertices $v,\,w$, we assign two {\em covering} edges
$e_{1},\,e_{2}$, where the edge $e_{i}$ connects $v_{i}$ with
$w_{i}$. We shall call the edges $e_{1}$ and $e_{2}$ {\em dual}.

We do the same with the {\em good} edges which do not belong to the
spanning tree $T$.

To each bad edge $e$, which connects the vertices $v$ and $w$, we
assign two {\em covering} edges $e_1,\,e_2$, the first edge connects
$v_{1}$ and $w_{2}$, the second edge connects $v_{2}$ and $w_{1}$.
The edges $e_1,\,e_2$ will be called {\em dual}.

Let $K^{2}$ be the constructed graph. The framing of a vertex of
$K^{2}$ is naturally induced from the graph $K$ which is covered by
$K^{2}$. One checks directly that the constructed graph $K^{2}$ and
the duality relation do not depend on the choice of the spanning
tree.

If $K$ is not connected, then the graph $K^{2}$ is constructed as
the split sum of the graphs $K_{i}^{2}$ that correspond to the
connected components $K_{i}$ of the graph $K$.

There is a natural involution on the graph $K^{2}$ that maps each
vertex to the dual vertex and each edge to the dual edge. Since the
duality is compatible with the framing relation at a vertex, there
is a natural definition of dual components of $K^{2}$.

The diagram $K'$ can be obtained from $K^{2}$ by removing one of the
two sets of components.

The following statement can be checked straightforwardly.
\begin{thm}\cite{ManSb}
The maps $f:K\to K', d:K\to K^{2}$ are well defined, i.e., when applying the Reidemeister move
to the diagram $K$, the diagrams  $K', K^{2}$ are operated on by combinations of Reidemeister moves.
\end{thm}

\begin{rk}
In the case if $K$ contains of a unique component, the map
$K\to K'$ is described as follows. The chord diagram $D(K')$ is obtained from $D(K)$ by deleting all odd chords.
\end{rk}

\begin{figure}
\centering\includegraphics[width=200pt]{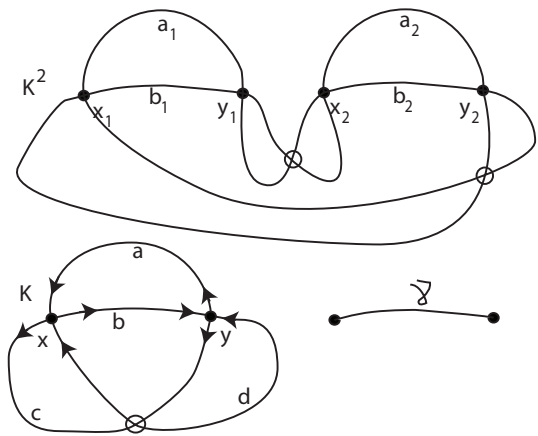}
\caption{The covering of a framed $4$-graph}
\label{covering}\end{figure}

An example of the covering $K^{2}$ over a framed $4$-graph $K$ is shown in Fig. \ref{covering}.

The graph $K$ shown in Fig. below, has two vertices  $x$ and $y$. For the spanning tree  $\gamma$, we shall take the tree consisting of one edge $b$ and two vertices $x$,$y$. The remaining edges will also connect $x$ to $y$. Herewith, the edge $a$ is {\em good}, and the edges
$c$ and $d$ are {\em bad}. In the upper part of Fig. \ref{covering}, we show the orientation at vertices $x$ and $y$, which agrees along the edge $b$.

This orientation is well defined on the edge $a$ but disagrees with itself along $c$ and $d$. The cycle $(b,a)$ is good since it {\em rotates} at one vertex and is not transverse anywhere, and both cycles $(b,c),(b,d)$ are bad (each of them is transverse at some vertex).

\section{Parities relative to a component of a free link}

The first known parity was the Gaussian parity. It was defined on
free knots and virtual knots and it is the only non-trivial parity
for free knots.

We shall see below that for two-component links one can define a
similar parity for one component (for crossings of the other
component the parity is not defined).

If we want to combine such a parity with the covering map, we
construct for a knot $K$ its covering $K^{2}=K_{1}\cup K_{2}$ and
get the parity on the crossings of the diagram $K_{1}$ relative to
the diagram $K_{2}$. Since the crossings of the diagram $K_{2}$ come
from the even crossings of the diagram $K$ the even crossing of $K$
split into two types. This construction was essentially described
in~\cite{ManSb2} where one discusses the ``refined parity'' of even
crossings in a diagram which contains both even and odd crossings.

Nonetheless, that refined parity is based on calculation of
intersections in the diagram.

In Fig.  \ref{cat2} we have demonstrated how to construct rather elaborated parities for
one component of a free knot by using sample subdiagrams.

It appears that by combining covering approach and the new parities
for one component of two-component link relative to the other
component, we can get very non-trivial parities for crossings of a
component $K_{1}$ in the link $K_{1}\cup K_{2}$ obtained as the
covering of a diagram $K$.

Such parities allow one immediately to construct and refine various
invariants of the initial diagram $K$, whose values depend on the
presence of one or another sample.

\subsection{Definition. The Simplest Example}

Let  $K\cup L$ be a two-component free link with evenly many mixed crossings.

Let us call a third Reidemeister move of the link $K\cup L$ {\em
special} if it involves one crossing of the component $K$ and two
mixed crossings (of the component $K$ with the component $L$).

\begin{figure}
  \centering\includegraphics[width=200pt]{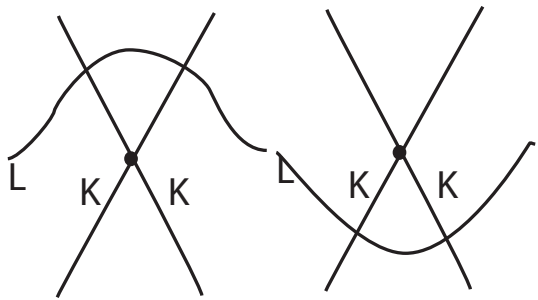}
  \caption{Mixed move does not change the parity}\label{KL}
 \end{figure}

 \begin{dfn}
Let $p$ be a rule which for some class of free links $K\cup L$
assigns a number $p_{v}$ equal $0$ or $1$ to each crossing $v$ of
the component $K$ . We call $p_{v}$ {\em parity for $K$ in $K\cup
L$}\index{Parity} if the following axioms hold:

a) usual parity axioms (see
Definition~\ref{parityDef}) for the Reidemeister
moves which involve only crossings of the component $K$;

b) Reidemeister moves which involve crossings of the component $L$
or mixed crossings do not change parity of the crossings of $K$ (see
Fig.~\ref{KL});

c) any special third Reidemeister move does not change the parity of
the crossing of the component $K$.
 \end{dfn}

For each crossing $v$ of the component $K$ we define the parity
$p_{L}(v)$ as the parity of number of mixed crossings in a half
$K_{v}$ (the choice of the half does not matter).

\begin{thm}
$p_{L}$ is a parity for $K$ in $K\cup L$.
\end{thm}

The statement follows from the straightforward check of Reidemeister
moves.

\begin{figure}
\centering\includegraphics[width=200pt]{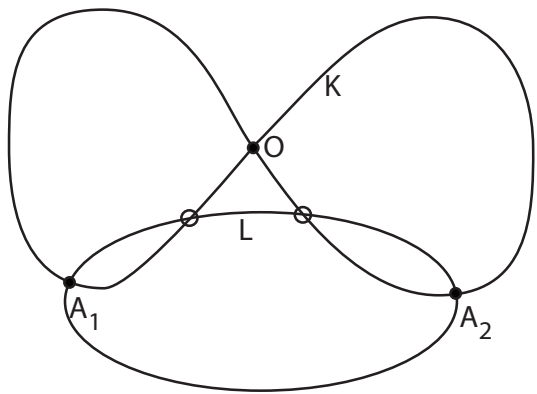}
\caption{A two-component link $K\cup L$} \label{kl2}
\end{figure}

\begin{example}
Let us consider the following two-component link $K\cup L$. It has one pure crossing $O$ belonging to the component $K$ and two mixed crossings $A_{1}$ and $A_{2}$; one of them belongs to one half of $K$ with respect to $O$, and the other one belongs to the other half of $K$ with respect to $O$.

According to the definition, the crossing $O$ is the only odd crossing
of the link diagram.

It follows from the definition that for any other diagram of the link $K\cup L$, the component $K$ also contains an odd crossing. These arguments yield the minimality of the diagram shown in Fig \ref{kl2}.
\end{example}

\subsection{Parities obtained from Turaev's cobracket}

The parity considered above for one component, being relative to the
other component, like the Gaussian parity, is defined by an
arithmetic calculation of crossings of some type. It appears that
one can construct parities which depend on the presence of some
``sublink'' in the given link. For this we apply some operation
(delta or bracket) to the component $L$ of a two-component link
$K\cup L$ and look how the result of the operation intersects the
component $K$.

Let us describe the parity construction in more detail.

The first parity will depend on some link $P\cup Q$, the presence of
the link as a resolution of the diagram means all the crossings of
the diagram are odd.

Let $P\cup Q$ be a free non-splittable two-component link whose
component $P$ and $Q$ are not trivial and are not equivalent to each
other.

Let ${\cal F}$ be the set of equivalence classes of non-splittable
free three-component links and let ${\cal M}={\Z}_{2}{\cal F}$ be
the corresponding free module.

We consider a map $f_{P,Q}$ which relates ordered two-component
links with an element of the module ${\cal M}$. This map is the
composition of two maps: $f_{P,Q}(K\cup L)=pr_{\cdot,P,Q}\circ
\Delta_{L}$.

Here $\Delta_{L}$ is the ordered map $\Delta$ which maps
two-component links to $\Z_{2}$--linear combinations of
three-component links. This map relates a link $K\cup L$ to the sum
$\sum_{s} K\cup L_{s,1}\cup L_{s,2}$, where the resolution $L\to
L_{s,1}\cup L_{s,2}$ is taken for all pure crossings $s$ of the
diagram $L$, and we do not consider splittable summands $K\cup
L_{s,1}\cup L_{s,2}$; as a result we get a linear combination of
three-component links with a distinguished component $K$.

By applying the projection $pr_{\cdot, P,Q}$, we keep only those
summands for which the two-component link obtained from $L$ is
equivalent to $P\cup Q$. Since $P$ and $Q$ are not equivalent the
components of each summand can be ordered so that the first
component is $K$, the second component is $P$ and the third
component is $Q$.

\begin{st}
The map $\Delta_{L}$ is well defined.
\end{st}
The statement follows from the direct check of Reidemeister moves.

Let us define the parity $p_{P,Q}$ for the knot $K$ relative to $L$
in the link $K\cup L$.

Let $K\cup L$ be a link with even number of mixed crossings. Each
pure crossing $v$ of the component $K$ of the link $K\cup L$
corresponds naturally to a crossing in every term $K\cup L_{1}\cup
L_{2}$ in the sum $f_{P,Q}(K,L)$. We denote this crossing as $v$ as
well. One can consider the two-component sublink $K\cup L_{1}$ and
the parity $p_{P}(v)$ relative to the component $L_{1}$ equivalent
to $P$. Summing these parities for all the summands of the sum
$f_{P,Q}(K,L)$, we get the number $p_{P,Q}(v)$.

\begin{st}
For links $K\cup L$ with even number of mixed crossings, the parity $p_{P,Q}(v)$ is a well-defined parity for $K$ relative to
$L$.
\end{st}

\begin{proof}
For all Reidemeister moves except for a second move applied to the
component $L$ the number of summands in $\Delta$ does not change,
one can identify all crossings of $K$ in every summand, and the
corresponding parities do not change.

A second move on the component $L$ leads to two new summands. Their
contributions into the parity for each crossing of the diagram $K$
annihilate, see Fig.~\ref{twocrossinginvariance}.

The figure shows two diagrams that differ with a second Reidemeister
move. Below the two contracting summands are shown.

\begin{figure}
\centering\includegraphics[width=200pt]{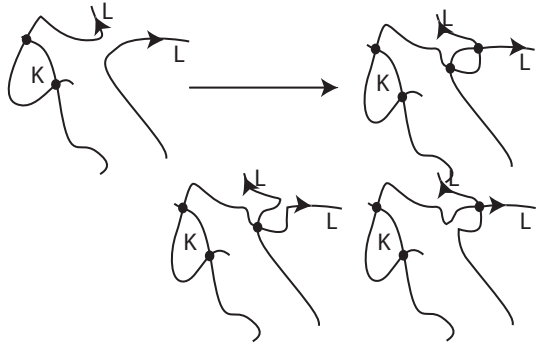}
\caption{Summands which annihilate by second Reidemeister move}
\label{twocrossinginvariance}
\end{figure}

\end{proof}

\subsection{Parities obtained from the parity bracket}

The idea of another parity construction is the following. Let $K\cup
L$ be a two-component link and let the component $L$ be odd and
irreducible. Consider $L$ as a four-valent graph and consider a
basis $\alpha_{1},\ldots,\alpha_{k}$ of the homology group
$H_{1}(L,\Z_{2})$ of the graph. Each element $\alpha_{i}$ is
presented as a cycle (or sum of cycles) in $L$. Assume that each of
these cycles has even number of intersections with the component
$K$.

Then we can call a crossing $v$ of the component $K$ of the diagram
$K\cup L$ to be {\em even} (resp., {\em
odd})\index{Crossing!even}\index{Crossing!odd} relative to
$\alpha_{i}$ if the intersection number of a component half at the
crossing $v$ with the cycle $\alpha_{i}$ is even (resp., odd).

Let us define a linear space $\gG''$ as the set of
$\mathbb{Z}_{2}$--linear combinations of the following objects. We
consider all framed 4-graphs with two ordered unicursal components
$K\cup L$ modulo the following equivalence relations:

1) all possible Reidemeister moves for the first component and all
possible mixed Reidemeister moves;

2) second Reidemeister moves on the component $L$.

\begin{rk}
Unlike $\gG'$, there is no obvious recognition algorithm for  $\gG''$.
\end{rk}

Let $K\cup L$ be the framed $4$-graph with two ordered components $K$ and $L$. Let us define the invariant $[K\cup L]_{L}$ according to the formula
 $$
[K\cup L]_{L}=\sum_{s_{{even}}}K\cup (L_{s_{even, 1 comp}})\in \gG'',
 $$
 where the sum is taken over all even splittings $s_{even}$ of the framed $4$-graph $L$
 having one component.

A straightforward Reidemeister move check leads us to the following
\begin{thm}
$[K\cup L]_{L}$ is an invariant of $2$-component links with ordered components.
\end{thm}

\section{Further discussion}

It turns out that this definition can be extended for a parity of
the component $K$ in $K\cup L$ by using $[K\cup L]_{L}$.

Let us now consider $[K,L]_{L}$.

It turns out that in some special cases, for each diagram ${\cal l'}$ of the link $L$, one can identify the homology classes $p_{1},\dots, p_{k}$ which correspond to the classes $\alpha_{1},\dots, \alpha_{k}$ of the diagram $L$ in such a way that  the following holds:
\begin{st}
Let us consider two-component links possessing the above properties. Then the maps $p_{1},\cdots, p_{k}$ are well defined parities for the knot $K$ with respect to $L$.
\end{st}

\end{document}